\documentclass[11pt]{article}

\usepackage[all]{xy}

\newtheorem{tm}{Theorem}[subsection]
\newtheorem{lm}[tm]{Lemma}
\newtheorem{pr}[tm]{Proposition}
\newtheorem{rmk}[tm]{Remark}

\newtheorem{??}[tm]{Question}
\newtheorem{defi}[tm]{Definition}

\newtheorem{ass}[tm]{Assumption}

\font\tenmsb=msbm10
\font\sevenmsb=msbm7
\font\fivemsb=msbm5
    
\newfam\msbfam
\textfont\msbfam=\tenmsb
\scriptfont\msbfam=\sevenmsb
\scriptscriptfont\msbfam=\fivemsb
\def\Bbb#1{{\fam\msbfam #1}}

\font\teneufm=eufm10
\font\seveneufm=eufm7
\font\fiveeufm=eufm5
\newfam\eufmfam
\textfont\eufmfam=\teneufm
\scriptfont\eufmfam=\seveneufm
\scriptscriptfont\eufmfam=\fiveeufm
\def\frak#1{{\fam\eufmfam\relax#1}}

\def\lorw{\longrightarrow}
\newcommand\n{\noindent}

\newcommand\ci{\cite}

\newcommand\rat{{\Bbb Q}}
\newcommand\comp{{\Bbb C}}

\newcommand\zed{{\Bbb Z}}

\newcommand\pn[1]{{\Bbb P}^{#1}}

\newcommand\blacksquare{{\hspace*{\fill} $\fbox{}$}}

\newcommand{\im}{ \hbox{\rm Im} }

\newcommand{\ke}{ \hbox{\rm Ker} }

\newcommand{\ptd}[1]{ \,^{\frak p}\!\tau_{ \leq {#1} } }

\newcommand{\ptu}[1]{ \,^{\frak p}\!\tau_{ \geq {#1} } }

\newcommand{\pe}{ {\cal P }  }

\newcommand{\pc}[2]{ \,^{\frak p}\!{\cal H}^{#1}({#2})   }

\newcommand{\ho}{\mbox{\rm Hom}}

\title{The perverse filtration  and \\  the Lefschetz hyperplane theorem}
\author{
Mark Andrea A.  de Cataldo\thanks{
Partially supported by N.S.F.}\, 
and Luca Migliorini\thanks{Partially supported by GNSAGA and PRIN 2007 project ``Spazi di moduli e teoria di Lie''}
}
\date{May 2008; Revised December 2008}

\begin{document}\maketitle

\begin{abstract}
We  describe
 the perverse filtration in cohomology 
using the Lefschetz hyperplane theorem.
\end{abstract}

\tableofcontents

\section{Introduction}
\label{intro}
In this paper we give a geometric description of the  middle perverse filtration on
the cohomology and on  the cohomology with compact supports 
of a complex   with constructible cohomology  sheaves
of abelian groups on a quasi projective variety. The description is in terms of restriction to generic hyperplane sections and it is somewhat unexpected, especially if one views
the constructions leading to perverse sheaves  as 
transcendental and hyperplane sections as more algebro-geometric.

The results of this paper are listed in $\S$\ref{results}, 
and hold for a  quasi projective variety.
For the sake of simplicity, 
we describe here the case 
of the cohomology of 
a $n$-dimensional affine variety $Y \subseteq {\Bbb A}^N$ with
coefficients in a  complex $K$.

The theory of $t$-structures endows the (hyper)cohomology groups
$H(Y,K)$ with a canonical filtration $P$, called the perverse filtration,
\[P^pH(Y,K)= \im\, \{H(Y, \ptd{-p} K) \lorw H(Y,K)\},\] which   is
 the abutment of the perverse spectral sequence.
Let \[Y_*= \{Y \supseteq Y_{-1} \supseteq  \ldots \supseteq 
Y_{-n}\}\] be a sequence of closed subvarieties; we call this data
a $n$-flag. Basic sheaf theory  
  endows $H(Y,K)$ with the so-called
flag filtration $F$, abutment of the  spectral sequence
$E_1^{p,q}= H^{p+q}(Y_{p}, Y_{p-1}, K_{|Y_{p}}) \Longrightarrow
H^*(Y,K)$. We have $F^p H(Y,K)= \ke \,\{ H(Y,K) \to H(Y_{p-1}, K_{|Y_{p-1}}) \}$.
For an arbitrary $n$-flag, the perverse and flag filtrations
are unrelated.

In terms of filtrations, the main result of this paper is that if
the $n$-flag is obtained using  $n$   hyperplane sections
in sufficiently general position, then
\begin{equation}
\label{1}
P^p H^j(Y,K) = F^{p+j} H^j(Y,K).
\end{equation}
More precisely, we construct a complex $R\Gamma (Y,K)$ endowed with two filtrations
$P$ and $F$ and we prove (Theorem \ref{tma}) that
there is a natural isomorphism in the filtered derived category
$DF(Ab)$ of abelian groups
\begin{equation}
\label{2}
(R\Gamma (Y,K), P) = (R\Gamma (Y, K), Dec (F)),
\end{equation}
where $Dec(F)$ is the shifted filtration associated with $F$. Then
(\ref{1}) follows from (\ref{2}).

Our methods seem to break down in the non quasi projective case and
also for other perversities.

The constructions and results  are amenable to mixed Hodge theory.
We offer the following   application:
let    $f:X \to Y$ be any map of algebraic varieties,
$Y$  be quasi projective and $C$ be a bounded
complex  with constructible cohomology sheaves on   $X.$ 
Then the perverse Leray spectral sequences
can be identified with suitable flag spectral sequences on $X.$ 
In the special case when $K=\rat_{X},$ we obtain
the following result due to  M. Saito:
the perverse spectral sequences for  $H(X, \rat)$ and $H_{c}(X,\rat),$
are spectral sequences of mixed Hodge structures.
Further Hodge-theoretic applications concerning the decomposition
theorem  are mentioned in
Remark \ref{marenostrum} and will appear in \ci{pfII}.

The isomorphism (\ref{2}) lifts to the bounded derived category
$D^b (\pe_Y)$
of perverse sheaves with rational coefficients.  This was the basis
of the proof of our results in an earlier version of this
paper. The present formulation,
which short-circuits $D^b(\pe_Y)$,  is based on the statement
of Proposition \ref{bfpr} which has been suggested to us by an anonymous referee.
We are deeply grateful for this suggestion.
The main point is that a suitable strengthening of the Lefschetz
hyperplane theorem yields  cohomological vanishings
for the bifiltered complex $(R\Gamma (Y, K), P,F)$ which
yield (\ref{2}). These vanishings are completely analogous
to the ones occurring for topological cell complexes and, for example, one can fit
the classical Leray spectral sequence of a fiber bundle in the framework
of this paper.

The initial inspiration for this work comes from Arapura's  paper \ci{ara}, 
which  deals with the standard filtration, versus the  perverse one. 
In this case, the flag has to be special: 
it is obtained by using high degree hypersurfaces
 containing the bad loci of the 
ordinary cohomology sheaves. The methods of this paper
are easily adapted  to that setting; see \ci{somm}.

The fact that the perverse filtration is  related to general hyperplane sections confirms, 
in our opinion, the more fundamental role played by perverse sheaves with respect to ordinary 
sheaves. The paper \ci{ara} has also directed us to the beautiful \ci{nori} and the seminal \ci{be}. 
The influence on this paper of the ideas contained in   \ci{be,nori} is hard to overestimate.

\bigskip
\n
{\bf Acknowledgments.}
It is a pleasure to thank D. Arapura,  A. Beilinson, M. Goresky, 
M. Levine, M. Nori 
for stimulating conversations.
The first author thanks 
the University of Bologna and I.A.S. Princeton
for their hospitality during the preparation of this paper.
The second author thanks the Centro di Ricerca Matematica E. De Giorgi in Pisa,
the I.C.T.P., Trieste, and I.A.S., Princeton
for their hospitality during the preparation of this paper.
Finally, we thank the referees for pointing out several inaccuracies
in an earlier version of the paper, and for very useful suggestions on 
how to make the paper more readable.

The first-named author dedicates this paper to Mikki$\&$Caterina.

\section{Notation}
A variety is a separated scheme of finite type over the field
of complex numbers $\comp$. A map of varieties
is a map of $\comp$-schemes.

The results
of this paper hold for sheaves of $R$-modules,
where $R$ is a commutative ring with identity
with finite global dimension, e.g. $R= \zed,$ $R$ a field, etc.
For the sake of exposition
 we work with $R=\zed$, i.e. with sheaves of abelian groups.
 
 The results of this paper hold, with routine
 adaptations of the proofs, in the case
 of varieties over an algebraically closed field
 and \'etale sheaves with the usual
 coefficients: $\zed/l^m\zed,$ $\zed_l,$ $\rat_l$,
 $\zed_l [E]$, $\rat_l [E]$ ($E \supseteq \rat_l$ a finite extension)
 and $\overline{\rat}_l.$

We do not discuss further these variants, except
 to mention that the issue
 of stratifications is  addressed
 in \ci{bbd}, $\S$2.2 and  $\S$6. 
 The term {\em stratification} refers
 to  an  algebraic Whitney stratification \ci{go-ma2}.
 Recall that any two stratifications admit a common refinement
 and that maps of varieties can be stratified. 
 
 Given a variety $Y$, there is the category
 $D_Y= D_Y(\zed)$ which is the full subcategory of the 
  derived category of the category  $Sh_Y$ of sheaves of abelian 
 groups whose  objects are the bounded  complexes with 
 constructible cohomology sheaves,
 i.e. bounded complexes
 $K$ whose cohomology sheaves ${\cal H}^i(K)$,
 restricted to  the strata of a 
 suitable stratification $\Sigma$ of $Y$,
  become locally constant with fiber a finitely generated
 abelian group.
 For a given $\Sigma,$ a complex
 with this property is called $\Sigma$-constructible.
 
 Given a stratification $\Sigma$ of $Y$, there are the
full subcategories
$D_Y^{\Sigma} \subseteq D_Y$   of complexes which are 
$\Sigma$-constructible.

Given a map  
$f:X \to Y$ of varieties, there are the  usual four  functors $(f^*,  Rf_*,Rf_!, f^!)$.  By abuse of  notation, we denote $Rf_*$ and $Rf_!$
simply by $f_*$ and $f_!$.
The four functors
preserve stratifications, i.e. 
if $f: (X, \Sigma') \to (Y, \Sigma)$ is stratified, 
then  $f_*, f_!: D^{\Sigma'}_X \to D^{\Sigma}_Y$
and 
$f^*, f^!:  D^{\Sigma}_Y \to D^{\Sigma'}_X.$

The abelian categories $Sh_Y$ and $Ab:=Sh_{pt}$ have enough injectives.
 The right derived functor of global sections  is denoted $R\Gamma (Y,-).$
Hypercohomology groups are denoted simply by $H(Y,K).$
Similarly, we have $R\Gamma_c(Y,-)$ and $H_c(Y,K).$
 
 We consider only  the {\em middle} perversity $t$-structure
 on $D_Y$ \ci{bbd}.  The truncation 
 functors are denoted 
 $\ptd{i} : D_Y \to \,  ^{\frak p}\!D_Y^{\leq i}$,
 $\ptu{j} : D_Y \to \,   ^{\frak p}\!D_Y^{\geq j}$,
 the heart $\pe_Y :=  {^{\frak p}\!}D_Y^{\leq 0}
 \cap {^{\frak p}\!}D_Y^{\geq 0}$
 is the abelian category of perverse sheaves on $Y$
 and we denote the perverse cohomology functors
 $^{\frak p}\!{\cal H}^i:= \ptd{i}\circ \ptu{i} [i]: D_Y \lorw \pe_Y$.
 
 \n
 The perverse $t$-structure is compatible with a fixed stratification,
 i.e. truncations preserve $\Sigma$-constructibility
 and we have  $^{\frak p}\!{\cal H}^i:= : D_Y^{\Sigma}
  \lorw \pe_Y^{\Sigma}$, etc.

 In this paper, the results we prove  in cohomology have a counterpart in
 cohomology with compact supports.
 If we employ field coefficients, then middle perversity
 is preserved by duality and  the results
 in cohomology are equivalent to  the ones in cohomology with compact
 supports by virtue  Poincar\'e-Verdier Duality.
 
 \n
   Due to the  integrality of the coefficients,
   middle perversity is 
 not preserved by duality; 
 see \ci{bbd}, $\S$3.3.
However, we can prove the results in cohomology and in 
 compactly supported cohomology using the same techniques.
 For expository reasons, we  often emphasize cohomology.

 Filtrations, on groups and complexes, are always {\em finite}, i.e. $F^iK =K$ for $i \ll 0$ and $F^jK =0$
 for $j \gg 0$, and {\em decreasing}, i.e. $F^iK \supseteq F^{i+1}K.$
 We say that $F$ has type $[a,b]$, for $a \leq b \in \zed,$ 
 if $Gr_F^p K \simeq 0$ for every $i \notin [a,b].$
  
  A standard reference for the filtered derived category
  $DF(A)$ of an abelian category is
  \ci{illusie}. 
  Useful complements can be found in
  \ci{bbd}, $\S$3 and  in \ci{be}, Appendix.
  We denote the filtered version of $D_Y$ by $D_YF$.
  The objects are filtered complexes $(K,F)$, with $K \in D_Y.$
  This is a full subcategory of $D^bF(Sh_Y).$
 
  We  denote a  ``canonical"  
  isomorphism with the symbol $``=."$

 \section{The perverse and flag spectral sequences}
 \label{tfss}
In this paper, we relate the perverse spectral sequences
with certain classical objects that we call flag spectral sequences.

In order to do so, we  exhibit these spectral sequences
as the ones associated with a collection of  filtered complexes of abelian groups.
These, in turn, arise by taking the global sections of
(a suitable injective model of)
the complex $K$ endowed with the    filtrations
$P,F, G$ and $\delta$ which we are about to define.

In this section,  starting with a variety $Y$ and a  complex
$K \in D_Y$, we construct the multi-filtered complex
$(K,P,F,G, \delta)$ and we list its relevant properties. 
By passing  to global sections, we identify the ensuing spectral sequences 
of filtered complexes with the perverse
and flag ones.

\subsection{$(K,P)$}
\label{kp}
The system of truncation maps
$\ldots \to  \ptd{-p} K \to \ptd{-p+1} K \to \ldots$
is isomorphic in $D_Y$ to a system
of inclusion  maps
$\ldots \to P^p K' \to P^{p-1} K' \to \ldots$, where
the filtered complex $(K',P)$ of injective type,
 i.e. all $Gr^p_P K'$, hence all $P^pK'$ and $K'$, have injective entries;
see \ci{bbd}, 3.1.2.7. 
The  filtered complex
$(K',P)$  is well-defined up to unique
isomorphism in the filtered  $D_YF$  by
virtue of \ci{bbd}, Proposition 3.1.4.(i) coupled with 
the second axiom, ``$\ho^{-1} =0$," of $t$-structures.  
We replace $K$ with $K'$ and  obtain  $(K,P).$
In particular, from now on, $K$ is injective.

\subsection{Flags}
\label{flags}
The smooth irreducible projective variety $F (N, n)$
of $n$-flags on $\pn{N}$
 parameterizes
 linear $n$-flags $ \frak F= \{ \Lambda_{-1} \subseteq
\ldots \subseteq  \Lambda_{-n}\}$, where $\Lambda_{-p}
\subseteq \pn{N}$ is a codimension
$p$ linear subspace.

A linear $n$-flag $\frak F$ on $\pn{N}$ is said to 
be {\em general} if it belongs to a suitable Zariski dense open
subset of the variety of flags $F(N,n).$ 
We say that a {\em pair of flags is general} if the same is true
for the pair with respect to $F(N,n)\times F(N,n).$
In this paper, this open set depends on the complex $K$
and on the fixed chosen embedding $Y \subseteq \pn{N}$.
We discuss this dependence in $\S$\ref{trbc}.

A  linear  $n$-flag $\frak F$ on $\pn{N}$ gives rise
to a {\em $n$-flag} on $Y\subseteq \pn{N}$, i.e. an increasing  sequence  of closed subvarieties of $Y$:
\begin{equation}
\label{nf}
Y_*= Y_* (\frak F ):  \quad   Y = Y_0 \supseteq Y_{-1} \supseteq 
\ldots \supseteq Y_{-n}, \qquad
Y_{p} : = \Lambda_p \cap Y.
\end{equation}
We set $Y_{-n-1}:= \emptyset$ and we have
 the  (resp., closed, open and locally closed)
embeddings:
\begin{equation}
\label{ijk}
i_p: Y_p \lorw Y, \qquad j_p : Y\setminus Y_{p-1} \lorw Y,
\quad k_p: Y_{p}\setminus Y_{p-1} \lorw Y.
\end{equation}

Let $h: Z \to Y$ be a locally closed embedding. There are the
exact  functor
$(-)_Z= h_! h^*(-)$, which preserves
$c$-softness, and the left exact functor $\Gamma_Z$,
which preserves injectivity and satisfies
$H_Z (Y, K)= H(Y, R\Gamma_Z K)$; see \ci{k-s}. If $h$ is closed, then 
$R\Gamma_Z = h_!h^!=h_*h^!;$  and, since $K$
is injective, $R\Gamma_Z K = \Gamma_Z K.$
If $Z' \subseteq Z$ is closed, then we have the distinguished
triangle $R\Gamma_{Z'} K \to R\Gamma_Z K \to R\Gamma_{Z-Z'} K
\stackrel{+}\to$  which, again by the injectivity
of $K$,  is the 
 triangle associated with  the exact sequence
$ 0 \to \Gamma_{Z'} K \to \Gamma_Z K \to \Gamma_{Z-Z'} K
\to 0.$

\subsection{$(K, F, G, \delta)$}
\label{kfgd}
We have constructed $(K,P)$ of injective type.
Let $Y \subseteq \pn{N}$ be an embedding of the
 quasi projective variety $Y$.
Let  ${\frak F}, {\frak F}'$
be two, possibly identical,  linear $n$-flags on $\pn{N}$
with associated flags, $Y_*$ and $Z_*$ on $Y.$
We denote the corresponding maps (\ref{ijk})
by $i',j',k'.$

We define  the three filtrations $F,G$ and $\delta$
 on $K$.
They are well-defined, up to unique isomorphism,
in the filtered $D_YF.$

\medskip
The  {\em flag  filtration
$F= F_{Y_*} = F_{Y_* (\frak F) }$},  of type $[-n,0]$,   is defined by setting  
$ F^p K:=K_{Y - Y_{p-1}}$:
\begin{equation}
\label{ttff}
0 \;\subseteq\;  K_{Y -Y_{-1}} \;\subseteq \; \ldots \; \subseteq \;
K_{Y-Y_{p-1} } \;  \subseteq  \; K_{Y- Y_{-n}} \;  \subseteq \; K.
\end{equation}

 The {\em  flag filtration} $G= G_{Z_*} = G_{Z_* (\frak F')}$,
 of type  $[0,n]$,
is defined by setting
$ G^p K:=\Gamma_{ Z_{-p}}K:$
\begin{equation}
\label{ttffg}
0 \;\subseteq  \;\Gamma_{Z_{-n}} K\; \subseteq \;  \ldots \; \subseteq \;
\Gamma_{Z_{-p} } K \; \subseteq  \; \Gamma_{Z_{-1}} \;
K  \; \subseteq   \; K.
\end{equation}

The {\em  flag filtration $\delta = \delta (F_{Y_*(\frak F)}
G_{Z_*({\frak F}')})$},  of type  $[-n,n]$,
is the {\em diagonal filtration}
defined by 
$\delta^p K = \sum_{i+j=p} F^i K\cap G^jK.$

\medskip
Note that one does not need injectivity to define the filtrations.
However, without this assumption, the resulting filtration
$G$ and $\delta$ would not be canonically defined in $D_YF.$
Moreover, injectivity yields $\Gamma_Z K = R\Gamma_Z K$ 
for every locally closed $Z \subseteq Y$, a fact we use
throughout
without further mention.

\subsection{The graded complexes associated with $(K,P,F,G, \delta)$}
\label{profpf}
Recall that
$
Gr^p_{\delta}  = \oplus_{i+j=p} Gr^i_F Gr^j_G 
$
and that Zassenhaus Lemma implies that
$Gr^i_F Gr^j_G = Gr_G^j Gr_F^i.$
Since the formation of $F$ is an exact functor,
the formation of $G$ is exact when applied to injective sheaves,
and
 injective  sheaves are $c$-soft, we have 
 \begin{equation}
 \label{bbb1}
 Gr^p_P K, \quad Gr^p_G K, \quad Gr^j_G Gr^p_P \qquad
 \mbox{are injective,}
 \end{equation}
 \begin{equation}
 \label{bbb2}
 Gr^p_F K, \quad 
 Gr^i_F Gr^p_P K, \quad
  Gr^i_F Gr^j_G  K, \quad 
 Gr^i_F Gr^j_G Gr^p_P K \quad 
Gr^p_{\delta} K, \quad
 Gr^a_\delta Gr^p_P K \qquad \mbox{are $c$-soft.}
 \end{equation}
 In particular, we have that $P^p K,$ $G^p K$ and $P^pK \cap
 G^jK$ are injective. We have 
 the  analogous $c$-softness statement for (\ref{bbb2}), 
 e.g.  the $F^i K \cap G^j  \cap P^p K$ are $c$-soft.
 By construction (\S\ref{kp}), the filtration $P$
splits in each degree, and the formation of $F$ and $G$ is compatible with direct sums. Hence, we have the following list of natural isomorphism:

\n
1)
$Gr^p_P K=\pc{-p}{K}[p]$, 

\n
2) $Gr^p_F K = K_{Y_p-Y_{p-1}}= {k_p}_! k_p^* K$, 

\n
3)
$Gr^p_G K = \Gamma_{Z_{-p} - Z_{-p-1}} K= 
{k'_{-p}}_* {k'}^!_{-p} K$. 

\n
4)
 $Gr^j_G Gr^i_F K=Gr^i_F Gr^j_G K =
( \Gamma_{Z_{-j}-Z_{-j-1}}   K)_{Y_{i}-Y_{i-1}}
 = (   {k'_{-j}}_* {k'}^!_{-j} K)_{Y_{i}-Y_{i-1}}$.

\n
5)  $Gr^p_{\delta} K = \oplus_{i+j=a} 
( \Gamma_{Z_{-j}-Z_{-j-1}}   K)_{Y_{i}-Y_{i-1}}$.

\n
6) $Gr^i_F Gr^p_P =(\pc{-p}{K}[p])_{Y_p-Y_{p-1}}$.

\n
7) $Gr^j_G Gr^p_P =
 \Gamma_{Z_{-j}-Z_{-j-1}}  (\pc{-p}{K}[p])=
{k'_{-j}}_* {k'}^!_{-j}   \pc{-p}{K}[p]$.

\n
8) $Gr^i_F Gr^j_G  Gr^p_P K = 
( \Gamma_{Z_{-j}-Z_{-j-1}}  (\pc{-p}{K}[p]))_{Y_{i}-Y_{i-1}}$.

\n
9) $Gr^a_{\delta} Gr^p_P K = \oplus_{i+j=a} 
( \Gamma_{Z_{-j}-Z_{-j-1}}  (\pc{-p}{K}[p]))_{Y_{i}-Y_{i-1}}$.

 \begin{rmk}
 \label{occissimo}{\rm 
 If the pair of flags is general, then (cf.  $\S$\ref{trbc}, or   \ci{be}, Complement to $\S$3)
$$
( \Gamma_{Z_{-j}-Z_{-j-1}}   K)_{Y_{i}-Y_{i-1}}
=
 \Gamma_{Z_{-j}-Z_{-j-1}}   (K_{Y_{i}-Y_{i-1}}).
$$
In general, the two  sides differ, for the rhs is zero on $Y_{i-1}.$}
 \end{rmk}

\subsection{$(R\Gamma (Y,K), P,F,\delta)$ and $(R\Gamma_c (Y,K),
P,G, \delta)$}
\label{iduecc}
Since $K$ is injective, we have $R\Gamma (Y, K) =
\Gamma (Y, K)$,
$R\Gamma_c (Y, K) =
 \Gamma_c (Y, K).$ We keep $``R"$ in the notation.
 
By applying the left exact functors $\Gamma$ and $\Gamma_c$,
we  obtain  the multi-filtered 
 complexes of abelian groups
 \begin{equation}
 \label{bb3}
 (R\Gamma (Y, K), P,F,\delta), \qquad (R\Gamma_c (Y, K),
 P, G, \delta),
 \end{equation}
 by setting, for example,
 $P^p R\Gamma (Y, K) := \Gamma (Y, P^p K)$, etc.

Since injective sheaves and $c$-soft sheaves are 
$\Gamma$ and $\Gamma_c$-injective, we have
\begin{equation}
\label{bft}
 R\Gamma (Y, 
\pc{-p}{K}[p]) = \Gamma (Y,  Gr_P^p K) =Gr_P^p \, \Gamma (Y,K),
\end{equation}
\begin{equation}
\label{bftc}
 R\Gamma_c (Y, 
\pc{-p}{K}[p]) = \Gamma_c (Y,  Gr_P^p K) =Gr_P^p \, \Gamma_c (Y,K),
\end{equation}
with analogous formul\ae $\,$ for the following
graded objects 
 \begin{equation}
 \label{btfcm}
 Gr^p_F, \, Gr^p_G,\,
Gr^a_F Gr^b_P,\, Gr^a_G Gr^b_P,\,
Gr^i_F Gr^j_G,\, Gr^p_{\delta},\,
Gr^i_F Gr^j_G Gr_P^b,\, Gr^p_{\delta}Gr_P^b.
\end{equation}

 \begin{rmk}
\label{rwij}
{\rm Though the formation of $F$ does not preserve injectivity,
one can always take filtered injective resolutions.
In that case, the resulting $F^p K$ is 
not exactly $K_{Y- Y_{p-1}}$, etc.,  but rather an injective resolution
of it. This would  allow to drop the  mention of
 $c$-softness. On the other hand, 
the $F$-construction is exact and formul\ae $\,$ like 
the ones in $\S$\ref{profpf}
are readily proved.}
\end{rmk}

\subsection{The perverse and flag spectral sequences}
\label{tpafss}
With the aid $\S$\ref{profpf},\ref{iduecc} it is immediate
to recognize the $E_1$-terms of the spectral sequences
associated with the filtered complexes 
$(R \Gamma (Y, K), P, F, \delta)$ and $( R \Gamma_c (Y, K), P,G, \delta).$

\begin{defi}
\label{f1f}
{\rm ({\bf Perverse spectral sequence and filtration})
The {\em perverse spectral sequence}  for  $H(Y, K)$ 
is  the spectral sequences
of the filtered complexes $(R\Gamma (Y, K), P)$:
\begin{equation}\label{psqw}
 E_1^{p,q} = H^{2p+q}(Y, \pc{-p}{K}) \Longrightarrow
 H^*(Y,K) 
 \end{equation}
and the   abutment is the  
 {\em perverse filtration} $P$ on $H^*(Y,K)$
 defined
by
\begin{equation}
\label{ppii}
P^p  H^*(Y,K) = \im\, \{H^*(Y, \ptd{-p}K) \lorw H^*(Y,K)\},
\end{equation}
Similarly, for $H_c(Y,K)$ using $(R\Gamma_c(Y,K), P).$
 }\end{defi}
 
 Let $f:  X \to Y$ be a map of algebraic varieties and $C \in D_X.$
 \begin{defi}
 \label{p12lss}
 {\rm
 The {\em perverse Leray spectral sequences}
 for  $H(X,C)$  ($H_c(X,C)$, resp.) are the corresponding
 perverse spectral sequences  on $Y$ for $K := f_*C$
 ($K:= f_! C$, resp.).
 }
 \end{defi}

Let $Y \subseteq \pn{N}$ be an embedding of the quasi projective
variety $Y$, $\frak F, {\frak F}'$  be two
 linear flags on $\pn{N}$ and $Y_*$
and $Z_*$ be the corresponding flags on $Y$.
\begin{defi}
\label{fls}
{\rm ({\bf Flag spectral sequence and filtration ($F$-version)})
The {\em $F$ flag spectral
sequence}  associated with  $Y_*$ 
is the spectral sequence associated with
the filtered complex
$(R\Gamma (Y,K),F)$:
\begin{equation}
\label{fgsq11}
E_1^{p,q} = H^{p+q}(Y, K_{Y_p-Y_{p-1}}) \Longrightarrow
H^*(Y,K)
\end{equation}
and its abutment is the  {\em flag filtration} $F= F_{Y_*}$
on $H^*(Y,K)$ defined
by
\begin{equation}
\label{fsq111}
F^p H^*(Y,K) = \ke \, \{ H^*(Y,K) \lorw
H^*(Y_{p-1}, K_{|Y_{p-1}} )  \}.
\end{equation}
}
\end{defi}


\begin{defi}
\label{flsg}
{\rm ({\bf Flag spectral sequence and filtration   ($G$-version)})
The {\em $G$ flag spectral
sequence}  associated with  $Z_*$ 
is the spectral sequence associated with
the filtered complex
$(R\Gamma_c (Y,K),G)$:
\begin{equation}
\label{fgsqg}
E_1^{p,q} = H^{p+q}_c(Y,  {k_{-p}}_* k_{-p}^! K) \Longrightarrow
H^*_c(Y,K)
\end{equation}
and its abutment is the  {\em flag filtration} $G= G_{Z_*}$
on $H^*_c(Y,K)$ defined
by
\begin{equation}
\label{fsqg}
G^p H^*_c(Y,K) = \im \, \{ \, H^*_c(Y,\Gamma_{Z_{-p}}K) \lorw
H^*_c(Y, K )  \, \}.
\end{equation}
}
\end{defi}

\begin{defi}
\label{flsd}
{\rm ({\bf Flag spectral sequence and filtration   ($\delta$-version)})
The {\em $\delta$ flag spectral
sequences}  associated with  $(Y_*, Z_*)$ 
are the spectral sequences associated with
the filtered complexes
$(R\Gamma  (Y,K),\delta)$ and
$(R\Gamma_c (Y,K),\delta)$.}
\end{defi}

\begin{rmk}
\label{senod}
{\rm 
We omit displaying these spectral sequences since,
due to Remark \ref{occissimo}, they do not have familiar
$E_1$-terms. If the pair of flags is 
general, or merely in good position wrt $\Sigma$ and each other
(cf. $\S$\ref{trbc}), then
we have equality in Remark
\ref{occissimo} 
and the $E_1$-terms take the following form
(we write $H_{c, Z}$ for $H_c \circ R \Gamma_Z$):
\begin{equation}
\label{bifss}
E_1^{p,q} = 
\bigoplus_{i+j=p} H^{p+q}_{Z_{-j} - Z_{-j-1}}
(Y, K_{Y_{i}
-Y_{i-1}}) \Longrightarrow 
H^* (Y, K)
\end{equation}
\begin{equation}
\label{bifssc}
E_1^{p,q} = 
\bigoplus_{i+j=p} H^{p+q}_{c,Z_{-j} - Z_{-j-1}}
(Y, K_{Y_{i}
-Y_{i-1}}) \Longrightarrow 
H^*_c (Y, K)
\end{equation}
and their abutments are the  {\em flag filtrations} $\delta
= \delta (Y_*, Z_*)$
on  $H(Y, K)$ and on $H^*_c(Y,K)$ defined
by 
\begin{equation}
\label{mnb}
\delta^p H^*(Y,K) = 
\im  \,\{\bigoplus_{i+j=p} \, H^*_{Z_{-j}}(Y, K_{Y-Y_i}) \lorw H^*(Y,K) \,
\},
\end{equation}
\begin{equation}
\label{mnbc}
\delta^p H^*_c(Y,K) = 
\im  \,\{\bigoplus_{i+j=p} \, H^*_{c,Z_{-j}}(Y, K_{Y-Y_{i}}) \lorw H^*(Y,K) \,
\}.
\end{equation}
}
\end{rmk}


 \subsection{The shifted filtration and spectral sequence}
 \label{ts1100}
We need the  notion and basic properties
(\ci{ho2})  of   the {\em shifted filtration}
for a filtered complex  $(L,F)$ in an abelian category. We make
the definition explicit in $Ab.$

The {\em shifted filtration} $Dec(F)$ on $L$  is:
$$
Dec(F)^p L^l := \{ x\in F^{p+l} K^l \, | \; dx \in F^{p+l+1} L^{l+1} \}.
$$
The {\em shifted spectral sequence of $(L,F)$}
is the one  for $(L, Dec(F))$ and we have
\begin{equation}
\label{id}
Dec(F)^p H^l(L) = F^{p+l} H^l(L),
\qquad
E_r^{p,q} (L, Dec (F) )  \; = \; E_{r+1}^{2p+q,-p} (L, F).
\end{equation}

\section{Results}
\label{results}
We prove results  for $Y$ quasi projective.
The  statements and the proofs are  more transparent 
when $Y$ is affine. We state and prove the results in the affine case
first.

The multi-filtered  complexes of abelian groups $(R\Gamma (Y, K), 
P, F, \delta)$ and $(R \Gamma_c (Y, K), P,G,  \delta)$,
which give rise to the spectral sequences and filtrations
we are interested in, are defined in $\S$\ref{tfss}.

\subsection{The results  over an affine base}
In this section $Y$ is {\em affine} of dimension $n$ and $K \in D_Y$.
Let $Y \subseteq \pn{N}$ be a fixed embedding
and ${\frak F}, {\frak F}'$ be a pair of linear
$n$-flags on $\pn{N}$.

\begin{tm}
\label{tma} {\rm{\bf (Perverse filtration on
cohomology for affine varieties)}}

\n
Let $\frak F$ be
  general. 
There is a natural isomorphism
in the filtered derived category $DF (Ab)$:
$$
(R\Gamma (Y, K), P) \; \simeq  \; 
(R\Gamma (Y, K), Dec(F))
$$
identifying the perverse spectral sequence
with the shifted  flag spectral sequence so that
$$
P^p H^l (Y,K) \; =\;     F^{p+l}     H^l(Y,K)  \; = \; 
\ke \, \{  \, H^l(Y,K) \lorw H^l(Y_{p+l-1}, K_{|Y_{p+l-1}}) \, \}.
$$
\end{tm}

\begin{tm}
\label{tmac} {\rm {\bf (Perverse filtration
on $H_c$ and affine varieties)}}

\n
Let $\frak F'$ be
 a general. 
 There is a natural isomorphism
in the filtered derived category $DF(Ab)$:
$$
(R\Gamma_c (Y, K), P, G) \; \simeq  \; 
(R\Gamma_c (Y, K), Dec(G), G)
$$
identifying the perverse spectral sequence
with the shifted  flag spectral sequence so that
$$
P^p H^l_c (Y,K) \; =\;     G^{p+l}     H^l_c(Y,K)  \; = \; 
\im \, \{  \, H^l_c( Y, R\Gamma_{Z_{-p-l}} K) \lorw 
H^l_c (Y, K) \, \}.
$$
\end{tm}

\medskip
In what follows, $f: X \to Y$ is an algebraic map, with $Y$ affine,
$C \in D_X$, and given a linear $n$-flag $\frak F$ on $\pn{N}$,
we denote by $X_*= f^{-1} Y_*$ the corresponding pre-image $n$-flag
on $X.$
\begin{tm}
\label{tmapl}{\rm{\bf (Perverse Leray and affine varieties)}}

\n
Let ${\frak F}$ be general.
The perverse Leray spectral sequence for $H(X, C)$
is  the corresponding   
shifted $X_*$  flag spectral sequence.
The analogous statement for $H_c(X,C)$ holds.
\end{tm}

\begin{rmk}
\label{dva}
{\rm
The $\delta$-variants of Theorems \ref{tma}, \ref{tmac}, \ref{tmapl}
for cohomology and for cohomology with compact supports,
hold for the  $\delta$ filtration as well, and with the same proof.
In this case one requires the pair of flags to be general.
}
\end{rmk}

\begin{rmk}
\label{dopo}
{\rm
Rather surprisingly,
the differentials of the perverse (Leray) spectral sequences
can be identified with the differentials
of a flag spectral sequence. In turn, these
are classical algebraic topology objects stemming
from a filtration by closed subsets, i.e. from the cohomology
sequences associated with  the triples $(Y_p, Y_{p-1}, Y_{p-2}).$
}
\end{rmk}

\subsection{The results over a   quasi projective  base}
\label{tritqpc}
In this section,  $Y$ is a quasi projective variety
of dimension $n$ and $K \in D_Y.$

There are several ways to state and prove
generalizations of the results in $\S$\ref{results}
to the quasi projective case. 
We thank an anonymous referee for suggesting this 
line of argument as an alternative
to  our original two arguments 
that used Jouanolou's trick (as in \ci{ara}), and  finite and affine   \v{C}ech coverings.
For an  approach via  Verdier's spectral objects see
 \ci{pfII}.

Let $Y$ be quasi projective and $Y\subseteq \pn{N}$ be a fixed 
{\em affine} embedding and $({\frak F}, {\frak F}')$ be a  pair of
linear $n$-flags on $\pn{N}$. The notion
of $\delta$ flag spectral sequence  is defined in 
Definition \ref{fls}; see also Remark \ref{senod}

\begin{tm}
\label{cqpp}
{\rm ({\bf Quasi projective case via two flags})}
Let the pair of flags be general.
There are natural  isomorphisms
in $DF (Ab)$
$$
(R\Gamma (Y, K), P ) \simeq 
(R\Gamma (Y, K), Dec(\delta)),
\quad
(R\Gamma_c (Y, K), P) \simeq 
(R\Gamma_c (Y, K), Dec(\delta))
$$ 
identifying the perverse and the shifted  $\delta$ flag
spectral sequence, inducing the identity
on the abutted  filtered spaces. 

\n
Moreover, if $f:X \to Y$ and $C \in D_X$ are given,
then  the perverse Leray  spectral sequences
coincide with the
shifted  $\delta$ flag  spectral sequences  associated with the
preimage flags on $X.$
\end{tm}

\section{Preparatory material}
\label{prep}
\subsection{Vanishing results}
\label{perora}
\begin{tm}
\label{cdav}
{\rm ({\bf Cohomological dimension of affine varieties})}

\n
Let $Y$
be affine and $Q \in \pe_{Y}$ be a perverse sheaf on $Y$. Then
$$
H^{r}(Y, Q) =0, \;  \forall r >0, \qquad
H^{r}_c(Y, Q) =0, \;  \forall r<0.
$$
\end{tm}
{\em Proof.} We give several references. 
The original proof of the first statement is due to Michael Artin \ci{sga4}, XIV
and is valid in the \'etale context.
\ci{go-ma2}, $\S$2.5:  proved  for intersection homology
with compact supports and with twisted coefficients
on a pure-dimensional variety; the reader can translate
the results in intersection cohomology and intersection cohomology
with compact supports  on a pure-dimensional variety;
a standard  devissage argument implies the result for a perverse sheaf
$Q$
on arbitrary varieties: $Q$ is a finite extension of intersection cohomology
complexes with twisted, not necessarily semisimple, coefficients
on the irreducible components.
\ci{bbd}, Th. 4.1.1:
the  case of $H$ is proved directly; the case of  $H_c$ is proved
 for field coefficients  by invoking  duality, however,
one can prove it directly and for arbitrary coefficients.
The textbook  \ci{k-s} proves it 
for Stein manifolds (see loc.cit. Theorem. 10.3.8); the general case follows
by  embedding $Y$
as a closed subset of an affine space $i:Y \to \comp^n$ and by applying the statement to the perverse sheaf $i_*Q$.
\blacksquare

\bigskip
Let $Y$ be quasi projective. Fix an {\em affine}
embedding $Y \subseteq \pn{N}$. 

Let $\Lambda, \Lambda' \subseteq \pn{N}$ be two  hyperplanes, 
 $H := Y \cap \Lambda \subseteq Y$ 
 and  $j: Y \setminus H \to Y \leftarrow H : i$  be the
 corresponding
open and closed immersions. 
Note that $j^!=j^*.$ Similarly, for $\Lambda'$.

 \begin{tm}
\label{swl}
{\rm ({\bf Strong Weak Lefschetz})}

\n
Let $Y$ be quasi projective and $Q \in \pe_{Y}.$
If $\Lambda$  is general, then
 $$
H^{r}(Y, j_{!}j^{!}Q) =0, \;\; \forall r < 0, \qquad
H^{r}_c(Y, j_{*}j^{*}Q) =0, \; \forall r > 0.
$$
Let $(\Lambda, \Lambda')$ be a general pair.
Then we have that $j_! j^* j'_* {j'}^! Q = j'_* {j'}^!   j_! j^* Q$
and
$$
H^r (Y,  j_! j^* j'_* {j'}^! Q)= H^r_c(Y,  j'_* {j'}^!   j_! j^* Q) = 0, 
\;\; \forall r \neq 0.
$$
\end{tm}
{\em Proof.} 
We give several references for the first statement.
 \ci{be}, Lemma 3.3; this proof is valid in the \'etale context.
The second statement is observed in \ci{be}, Complement to $\S$3.
 \ci{go-ma2}, $\S$2.5. 
  M. Goresky has informed
us that P. Deligne has  also proved this result (unpublished). 
\blacksquare

\medskip
We  include a sketch of the proof of this result, following 
\ci{be}, in  $\S$\ref{trbc},
where we   also complement
the arguments in \ci{be}  that we need
in the sequel of the  paper.

\begin{rmk}
\label{falt}
{\rm
Since $j^! = j^*$, we may reformulate  the first statement
of Theorem \ref{swl}
as follows
$$
H^{r}(Y, Q_{Y-H}) =0, \;\; \forall r < 0, \qquad
H^{r}_c(Y, R\Gamma_{Y-H}Q) =0, \; \forall r > 0
$$
and similarly for the second one.
Moreover, 
by Theorem \ref{cdav},
if $Y$ is affine, then the vanishing results hold for every $r \neq 0$.
}
\end{rmk}

\begin{rmk}
\label{essaff}
{\rm 
It is essential that the embedding $Y \subseteq
\pn{N}$ be  affine. For example, the conclusion does not hold
in the case when $Y = {\Bbb A}^2 \setminus \{0\} \subseteq \pn{2}$
and $Q= \zed_Y [2].$
}
\end{rmk}

\subsection{Transversality, base change and choosing good flags}
\label{trbc}
In this section we highlight the role of  transversality 
in the proof of Theorem \ref{swl}. In fact,
transversality implies several base change equalities
which we use throughout the paper
in order to prove the vanishing results in
 Theorem \ref{swl},  its iteration
Lemma \ref{cellu},  its  ``two-flag-extension" (\ref{dova})
and to observe
 (\ref{bctcms}). While the vanishing results are used to realize 
condition (\ref{sta}), which is the key to the main results
of this paper,    the base change equality (\ref{bctcms}) is used to
reduce the results
for the perverse spectral sequences Theorems
\ref{tmapl}, \ref{cqpp} 
wrt a map $X \to Y$,
to   analogous results for perverse spectral sequences on $Y$.
 
These  base change properties hold generically by virtue
of the generic base change theorem \ci{sga4m}, and this is enough
for the purpose of this paper. On the other hand,
it is possible   to 
pinpoint the conditions one needs to impose on flags;
see Definition \ref{basdef} and Remark \ref{compno}.

\bigskip
 Let $Y \subseteq \pn{N}$ be an affine embedding of the quasi
projective variety $Y$ and 
$Y \subseteq \overline{Y}$
be  the resulting projective completion. There is a natural
decomposition  into locally closed subsets
$
\pn{N}\; =\; (\pn{N} \setminus \overline{Y}) \; \coprod\;
(\overline{Y}\setminus Y)\; \coprod \;Y .
$
Let $K \in D_Y.$
\begin{defi}
\label{adpted}{\rm
({\bf Stratifications adapted to the complex and to the embedding})
We say that
a stratification $\Sigma$ of $\pn{N}$ is {\em adapted to
the embedding $Y \subseteq \pn{N}$}  if  $Y,$  
$\overline{Y} \setminus Y,$ hence $\overline{Y}$, and
 $\pn{N} \setminus \overline{Y}$
are union of strata, and $\Sigma$ induces by restriction
stratifications on $\pn{N}, \pn{N} \setminus
\overline{Y}, \overline{Y}, Y, \overline{Y} \setminus Y$
 with respect to which
all possible  inclusions among these varieties are stratified maps. 
We denote these induced stratifications by $\Sigma_Y,$
etc.

\n
We say that $\Sigma$ is {\em adapted to $K$} if $K$ is $\Sigma_Y$-constructible.
}\end{defi}

\begin{rmk}
\label{rase}{\rm
Since maps of varieties can be stratified and a finite collection
of stratifications admits a common refinement, 
stratifications  which are adapted to the the complex and the embedding
exist.
}\end{rmk}

Let  $\Sigma$ be  a stratification of $\pn{N}$ adapted to $K$ and to the embedding
$Y \subseteq \pn{N}$.
Let  ${\Lambda} \subseteq {\Bbb P}$ be a hyperplane,
$H : = \Lambda \cap  Y$ and $ \overline{H}= 
\Lambda \cap \overline{Y}$.  Set
 $ \overline{U}:= (\overline{Y} \setminus
\overline{H})$ and $U := (Y \setminus H).$  Consider the cartesian  diagram
\begin{equation}
\label{su}
\xymatrix{
 H \ar[r]^i \ar[d]^J & Y \ar[d]^J & U  \ar[d]^J \ar[l]_j \\
\overline{H} \ar[r]^i & \overline{Y} & \overline{U}  \ar[l]_j.
}
\end{equation}
We address the following question:
when is the natural map 
\begin{equation}
\label{ooii}
J_! j_* j^* K \to j_* J_! j^*K
\end{equation}
an isomorphism? 
In general the  two differ on $\overline{H} \cap (\overline{Y} \setminus Y).$ By the octahedron axiom,  the map (\ref{ooii}) is an isomorphism
iff the natural base change map  $J_* i^*  K \to i^* J_* K$
is an isomorphism. This latter condition is met if $\Lambda$
is general (\ci{be}, Lemma 3.3). In fact it is sufficient
that $\Lambda$ meets 
transversally the strata in $\Sigma_{\overline{Y} -Y}$.
This is a condition on the stratification, not on $K$.
It follows that the analogous map
$j_! J_* J^! K \to J_* j_! J^! K$ is also an isomorphism
under the same conditions.

\medskip
\n
{\bf
Proof of Theorem \ref{swl}} (see \ci{be}).
 We prove the  first statement for cohomology.
The point is  that a general linear section produces the  isomorphism (\ref{ooii})
and this identifies the  cohomology groups
in question
with   compactly supported cohomology groups on affine varieties where one  uses Theorem \ref{cdav}.
Note that since the maps of type $j$ and $J$ are affine,
all the complexes appearing below are perverse. 
We have the following chain of equalities:
$$
H^r (Y, j_! j^* Q)= H^r (\overline{Y},J_* j_! j^* Q) =
H^r(\overline{Y}, j_! J_* j^* Q)   = H^r_c(\overline{Y}, j_! J_* j^! Q)
=   H^r_c(\overline{U},  J_* j^! Q)
$$
and, since $\overline{U}$ is affine and  $ J_* j^! Q$ is perverse,
the last group
is zero for $r <0$ and the first statement for cohomology follows.
The one for compactly supported cohomology
is proved in a similar way.

\n
In order to prove the second statement, we consider the Cartesian diagram
\begin{equation}
\label{sud}
\xymatrix{
 U\cap U' \ar[r]^{j'} \ar[d]^j & U \ar[d]^j  \\
U' \ar[r]^{j'} &  Y.
}
\end{equation}
Since the embedding $Y\subseteq \pn{N}$ is affine,
the open sets $U, U'$ and $U\cap U'$ are affine. Note that
this fails if the embedding is not affine.
We have that $j_!, j_*, j^!=j^*$ are all $t$-exact and preserve
perverse sheaves. The same is true for $j'$.

\n
The equality  $j_! j^* j'_* {j'}^! Q = j'_* {j'}^!   j_! j^* Q$ is proved using 
base change considerations  similar to the ones
we have made for (\ref{ooii}). 

\n
We prove the vanishing in cohomology.
The case of cohomology with compact supports
 is proved in a similar way. The case $r<0$ is
covered by the first statement. We need suitable
``reciprocal" transversality conditions which are the obvious
generalization of the ones mentioned when discussing
(\ref{ooii}). We leave the formulation of these conditions to the reader.
 It will suffice to say that they are met by a general pair $(\Lambda, 
 \Lambda').$
 The case $r>0$ follows from Theorem
\ref{cdav} applied to the affine $U'$:
$H(Y,  j'_* {j'}^!   j_! j^* Q) = H(U',   {j'}^!   j_! j^* Q).$
\blacksquare

\begin{rmk}
\label{mm0}
{\rm Let $Q \in \pe_Y$ and
 $\Sigma$ be a stratification of $\pn{N}$ adapted to $Y \subseteq \pn{N}$ and such that $Q \in \pe_Y^{\Sigma_Y}.$
An inspection of the proof of Theorem \ref{swl} reveals that it is sufficient
to choose $\Lambda$ so that it meets transversally
the strata in $\overline{Y} \setminus Y.$ It is not relevant
how $H$ meets the strata in $\Sigma_Y.$ A similar remark
holds in the case of a pair of hyperplanes.

}\end{rmk}

We now introduce a kind of transversality notion
that is sufficient for the purpose of this paper.
Let 
 $\Sigma$  be as above. 
 
\begin{defi}
\label{basdef}
{\rm
({\bf Flag in good position wrt $\Sigma$})

\n
A linear $n$-flag $\frak F= \Lambda_0 \subseteq \Lambda_{-1} \subseteq
\ldots \subseteq \Lambda_{-n}$ on $\pn{N}$ 
 is in {\em good position with respect
to $\Sigma$} if it is subject to the following inductively
defined conditions:

\n
1) $\Lambda_{-1}$ meets all the strata of $\Sigma_0:= \Sigma$ transversally;

\n
let $\Sigma_{-1}$ be a refinement of $\Sigma$ such that
its restriction to $\Lambda_{-1} \simeq \pn{N-1}$
is adapted to the embedding $Y_{-1} \subseteq \Lambda_{-1}$;

\n
2) $\Lambda_{-2}$ meets all the strata of $\Sigma_{-1}$ transversally;

\n
we iterate these conditions and constructions and introduce
$\Sigma_{-2},$ $\Lambda_{-3}$, \ldots, $\Sigma_{-n+1},$
$\Lambda_{-n}$  and we require
that, for every $i=1, \ldots, n$;

\n
i) $\Lambda_{-i}$ meets all the strata of $\Sigma_{-i+1}$ transversally.

\n
We define $\Sigma':= \Sigma_{-n}.$
}\end{defi}

\begin{rmk}
\label{fac}
{\rm
By the Bertini theorem, it is clear that  a general linear  $n$-flag $\frak F$  is in good position wrt
a fixed $\Sigma.$ Of course, ``general" depends on $\Sigma.$
Note also that if  $\frak F$ is in good position, 
then  $Y_p$ has pure codimension $p$ in $Y$.}
\end{rmk}

\begin{rmk}
\label{compno}
{\rm
There is the companion notion of a pair of linear $n$-flags 
$({\frak F}, {\frak F}')$ being {\em in good position
wrt to $\Sigma$ and  each other}. 
We leave the task of writing down the precise formulation
to the reader.  The notion is again inductive and proceeds,
also by imposing mutual transversality,  in
the following order: $\Lambda_{-1}, \Lambda_{-1}',
\Lambda_{-2}, \Lambda_{-2}'$, etc. It suffices to say
that   a general pair of  flags will do.
}
\end{rmk}

\smallskip
The proof of Theorem \ref{swl} works well inductively with the elements
of a linear flag $\frak F$ on $\pn{N}$
 in good position wrt to the embedding and to the perverse sheaf
$Q$.
We use this fact in the proof of Lemma \ref{cellu}.
Similarly, 
This kind of argument works well with a pair of flags in good position
with respect to the $Q$,  the embedding and each other.
In particular, it works for a general
pair of flags. In these cases,  we have the  equality
$
\Gamma_{Z_{-j} - Z_{-j-1}} (Q_{Y_i -Y_{i-1}}) =
(\Gamma_{Z_{-j} - Z_{-j-1}}Q)_{Y_i -Y_{i-1}}  
$ which follows from a repeated use
of the equality $j_! j^* j'_* {j'}^! Q = j'_* {j'}^!   j_! j^* Q$
 of Theorem \ref{swl}.
 By transversality, the shift $[i+j]$ of these  complexes 
are perverse.
This allows to apply the vanishing
results of Theorem \ref{swl} and deduce,
for general pairs of flags on the {\em quasi projective} 
variety $Y$ that
\begin{equation}
\label{dova}
H^r(Y,  \Gamma_{Z_{-j} - Z_{-j-1}} (Q_{Y_i -Y_{i-1}}))
= 
H^r_c(Y,  \Gamma_{Z_{-j} - Z_{-j-1}} (Q_{Y_i -Y_{i-1}}))
=0, \qquad \forall \, r \neq 0.
\end{equation}

\medskip
Let $f: X \to Y$ be a map of varieties.
The diagram (\ref{su}) induces
the cartesian diagram:
\begin{equation}
\label{alceste}
\xymatrix{
X_H \ar[r]^i \ar[d]^f  & X \ar[d]^f & X_U \ar[d]^f \ar[l]_j \\
H \ar[r]^i & Y & U \ar[l]_j.
}
\end{equation}
The previous  base change discussion  implies, for $\Lambda$ meeting all strata
of $\Sigma$ transversally, that 
\begin{equation}
\label{bctcms}
f_* j_! j^*C =  j_! j^* f_*C, \qquad 
f_! j_* j^! C =  j_* j^!f_! C.
\end{equation}
Similar base change equations hold for a
linear flag $\frak F$ on $\pn{N}$  in good position with respect to 
$f_!C$  and  to the embedding $Y \subseteq \pn{N}$
(e.g. general) 
and also for a pair of flags in good position
with respect to $f_! C$, the embedding and each other (e.g. a
 general pair).

\subsection{Two short exact sequences}
\label{ases}
\begin{lm}
\label{brs}
Let   $Y\subseteq \pn{N}$ be quasi projective and $Q \in \pe_Y.$
If $\Lambda \subseteq \pn{N}$ is  a general linear section, then
there are natural exact sequences in $\pe_Y$:
\begin{equation}
\label{brs1}
0 \lorw i_* i^*Q [-1] \lorw j_! j^! Q \lorw Q \lorw 0,
\end{equation}
\begin{equation}
\label{brs2}
0 \lorw  Q    \lorw j_* j^* Q \lorw   i_! i^!Q [1]\lorw 0.
\end{equation}
\end{lm}
{\em Proof.}
There are  the  distinguished triangles in $D_Y$:
$$
j_! j^! Q  \lorw Q \lorw  i_* i^* Q \stackrel{+}\lorw,
\qquad \qquad
i_! i^! Q  \lorw Q \lorw  j_* j^* Q \stackrel{+}\lorw.
$$
Since $j$ is affine,  $j_!$ and $j_*$ are $t$-exact and $j_! j^! Q$
and $  j_* j^* Q$ are perverse.
We choose, $\Lambda$ so that is it transverse to the strata
of  a stratification for $Q$. It follows that
$i_*i^*Q[-1] =i_! i^! Q [1]$ is  perverse. Each conclusion
follows from the long exact sequence
of perverse cohomology
of the corresponding   distinguished triangles.
\blacksquare


\subsection{The forget the filtration map}
\label{tftfm}
Let $A$ be an abelian category.
\ci{bbd}, Proposition 3.1.4.(i)   is a sufficient
condition for the natural forget-the-filtration map
$\ho_{DF(A)} \to \ho_{D(A)}$ to be an isomorphism.
We need 
the bifiltered counterpart of this sufficient condition.

The objects
$(L,F,G)$ of the bifiltered derived category $DF_2(A)$ are  complexes
$L$ endowed with two filtrations. The homotopies
must respect both filtrations and one inverts bifiltered quasi isomorphisms,
i.e. (homotopy classes of) maps inducing  quasi isomorphisms on the bigraded objects $Gr^a_F Gr^b_G.$ 
 It is a routine matter
to adapt Illusie's treatment of  $DF(A)$ to the bifiltered setting
and then to adapt  the proof of  \ci{bbd}, Proposition 3.1.4.(i)
to yield a proof of 
\begin{pr}
\label{sv}
Assume that
 $A$ has enough injectives and that
  $(L,F,G), (M,F,G) \in D^+F_2(A)$ are such that 
$$
\ho_{DFA}^n ( (Gr_G^iL[-i],F),  (Gr_G^jM[-j],F)) =0, \qquad 
\forall  n<0, \quad \forall i<j.
$$
The ``forget the second  filtration" map is an  isomorphism:
$$
\ho_{DF_2A}((L,F,G),(M,F,G))  \stackrel{\simeq}\lorw \ho_{DFA}(
 (L,F),(M,F)).
$$
\end{pr}

\subsection{The canonical lift of a $t$-structure}
\label{canlift}
The following is a mere special case of
\ci{be}, Appendix.

Let $A$ be an abelian category.
The  derived category $D(A)$ admits the standard $t$-structure,
i.e. usual truncation.
The filtered derived category $DF(A)$  admits a canonical $t$-structure which
lifts (in a suitable sense which we do not need here)
 the  given one on $D(A)$.
 This canonical  $t$-structure on $DF(A)$  is 
described as follows.
There are the two full subcategories
$$
DF(A)^{\leq 0}: =  \{ (L,F) \, | \; Gr^i_FL  \in D(A)^{\leq i} \},
\quad
DF(A)^{\geq 0}: =  \{ (L,F) \, | \; Gr^i_FL  \in D(A)^{\geq i} \}.
$$
The heart  $DF(A)^{\leq 0} \cap DF(A)^{\geq 0}$ is   
$$
DF_{\beta} (A) = \{ (L,F) \, | \;  Gr^i_F L [i]  \in A \},
$$
where $\beta$ is for b\^{e}te (see \ci{bbd}, 3.1.7).
The reader can verify the second axiom  of $t$-structure,
i.e.
$\ho^{-1} (DFA^{\leq 0}, DFA^{\geq 1})=0$, by 
a simple induction on the length of the filtrations,
and the third axiom, i.e. the existence of the truncation triangles,
by simple  induction on the length of the filtration  coupled
with the use of Verdier's  ``Lemma of nine" (see \ci{bbd}, Proposition
1.1.11).

\subsection{The key lemma on bifiltered complexes}
\label{lobc}
Let $A$ be an abelian category
and $(L, P,F)$ be a bifiltered
complex, i.e. an object in the bifiltered derived category 
$DF_2(A)$. Recall  the  existence of 
the shifted filtration $Dec(F)$ associated with 
$(L,F)$.

The purpose of this section is to prove
the following  result the formulation of which has been suggested
to us by an anonymous referee. This result is key
to the approach presented in this paper. 
\begin{pr}
\label{bfpr}
Let $(L,P,F)$ be a bifiltered complex  be such that
\begin{equation}
\label{sta}
H^{r}  (Gr^a_F Gr^b_P L) =0, \;\; \forall r \neq a-b.
\end{equation} 
Assume that $L$ is bounded below and that $A$ has enough injectives.

\n
There is a natural isomorphism in the filtered derived category $DF(A)$
$$
(L, P) \simeq (L, Dec(F) )
$$
that induces the identity on $L$ and thus identifies
$(H^*(L), P) = (H^*(L),Dec(F)).$

\n
In particular, there is a natural isomorphism between the spectral sequences
associated with  $(L,P)$ and $(L,Dec(F))$ inducing the identity on the abutments.
\end{pr}

In order to prove Proposition \ref{bfpr},
we need the following two lemmata.

\begin{lm}
\label{lmlu}
Let $(L,F)$ be any filtered complex.
Then  the bifiltered complex
 $(K,F,Dec(F)) $ satisfies (\ref{sta}).
\end{lm}
{\em Proof.} This is a formal  routine verification.
\blacksquare

\begin{lm}
\label{giso}
Let things be as in Proposition \ref{bfpr}.
The natural map 
$$
\ho_{DF_2(A)} ( (L,P, F), (L,Dec(F), F)) \lorw
\ho_{DF(A)}((L,F), (L,F))
$$
induced by forgetting the first filtration is an isomorphism.

\n
The same is true with the roles of  the filtrations $P$ and $Dec(F)$ switched.
\end{lm}
{\em Proof.} Endow $D(A)$ with the standard $t$-structure
(i.e. usual truncation). Endow $DF(A)$ with the canonical
lift of this $t$-structure (see $\S$\ref{canlift}).

\n
The hypothesis (\ref{sta}) implies that, for every $b\in \zed,$
$(Gr_P^bL[-b], F)    \in
DF_{\beta} (A)$, i.e. it is in the heart
of the canonical lift of the standard  $t$-structure to
$DF(A).$
Similarly, Lemma \ref{lmlu} implies that, for every $b\in \zed,$
$(Gr_{Dec(F)}^bL[-b], F) \in 
DF_{\beta} (A)$.

\n
The  hypotheses of  Proposition
\ref{sv} are met: in fact they are met
  for every $i,j$,
due to the second axiom of $t$-structure.
The first statement follows.

\n
If we switch $P$ and $Dec(F),$ then  the hypotheses
of Proposition \ref{sv}
are still met, for the same reason, and the second statement follows.
\blacksquare

\bigskip
\n
{\bf Proof of Proposition \ref{bfpr}.} 

\n
By Lemma $\ref{giso}$, the identity on $(L,F)$ admits 
natural lifts
$$\iota_P  \in \ho_{DF_2A} ( (L,P, F), (L,Dec(F), F)), \quad
\iota_{Dec(F)} \in  \ho_{DF_2A} ( (L,Dec(F), F), (L,P, F))$$
which are inverse to each other and hence isomorphisms.

\n
By forgetting the second filtration, we obtain
a pair of maps in $\ho_{DFA} ( (L,P), (L,Dec(F))$
$\ho_{DFA} ( (L,Dec(F)), (L, P)$ which are inverse to each other.

\n
By forgetting both filtrations,  both maps yield the identity
on $L$.
\blacksquare

\begin{rmk}
\label{arci}
{\rm
The results of this section hold if we replace $Dec(F)$ with any filtration
$P'$ satisfying (\ref{sta}).
}
\end{rmk}

\section{Proof of the results}
In this section, we prove the main results of this paper and
we make a connection with Beilinson's equivalence theorem \ci{be}.

\subsection{Verifying the vanishing 
(\ref{sta}) for general flags}
\label{vvsfgf}
Recall the set-up: $Y$ is quasi projective
of dimension $n$,  $K \in D_Y$, 
$Y \subseteq \pn{N}$ is an affine embedding,
${\frak F}, {\frak F}'$ is a  pair of linear $n$-flags
on $\pn{N}.$
 We have  the bounded  multi-filtered complexes
 in of abelian groups
$ (R\Gamma (Y,K), P,F, \delta) ,  \, (R\Gamma_c (Y,K), P,G,\delta)$
 obtained using suitably acyclic resolutions.
  The perverse spectral sequences are the spectral sequence
for the filtration $P$, the flag spectral sequences
are  the ones for the filtrations $F, G, \delta$. Similarly,
for the perverse Leray spectral sequences.
If the flags  are arbitrary, then the perverse
and the flag   spectral sequences
seem unrelated. 

Let $\Sigma$ be a stratification of $\pn{N}$ adapted to 
$K$ and to the embedding $Y \subseteq \pn{N}.$
The proof of  Theorem \ref{tma}    consists
of showing that 
if the flag $\frak F$ is in   good position wrt $\Sigma$
(see Definition \ref{basdef}),
then the vanishing conditions  (\ref{sta}) hold for
the bifiltered complexes $(R\Gamma (Y, K), P, F))$ by virtue of a repeated application
of  the strong weak Lefschetz theorem \ref {swl}, so that
Proposition \ref{bfpr} applies and
there is a natural identification of filtered complexes $(R\Gamma (Y, K), P) = 
(R\Gamma (Y, K), Dec (F))$ and of the ensuing spectral sequences.
The other results are proved in a similar way.

The key to the proof is
Lemma \ref{cellu} which is  suggested by 
a construction  due to Beilinson
\ci{be}, Lemma 3.3 and Complement to $\S$3.
It is a technique to find resolutions of perverse sheaves 
on  varieties by using suitably transverse flags. 
The entries of the resolutions
satisfy  strong vanishing conditions and  realize 
the wanted condition (\ref{sta}). There are three
versions, left, right and bi-sided resolutions. 
 
 The resolutions are complexes obtained through the following general
 construction. Let $Q \in \pe_Y^{\Sigma_Y}$, where $\Sigma_Y$
 is the trace of $\Sigma$ on $Y$.
 The connecting maps associated with
 the short exact sequences  $0\to Gr_F^{*+1} Q \to
 F^{*}Q/F^{*+2} Q \to Gr^* Q \to 0$ give rise to a sequence 
 of maps in $D_Y$ 
 \begin{equation}
 \label{cfses}
 Gr^{-n}_F Q \stackrel{d}\lorw \ldots  \stackrel{d}\lorw Gr^0_F Q.
\end{equation}
with $d^2=0.$ We call this a complex in $D_Y$.
The same is true for the $G$ filtration:
$Gr^0_G Q \to \ldots \to Gr^n_G Q$ is a complex in $D_Y$.
The bigraded objects $Gr_F Gr_G$ give rise to a double complex
with associated single complex $Gr^{-n}_{\delta} Q \to 
\ldots  \to Gr^n_{\delta} Q$ in $D_Y.$
The transversality assumptions on the flags  ensure
that these 
are complexes of perverse sheaves resolving $Q$, 
that they 
are  suitably acyclic and  that their formation
is an exact functor. More precisely, we have the following

\begin{lm}
\label{cellu}
{\rm ({\bf Acyclic resolutions of perverse sheaves})}
Let $Y$ be quasi projective, 
$\Sigma$ be a  stratification  adapted to the  affine embedding
 $Y \subseteq \pn{N}$,  
 $Q \in \pe_Y^{\Sigma_Y}$ be a $\Sigma_Y$-constructible
 perverse sheaf on $Y$.
 
 \n
 Let $\frak F$ be a  linear  $n$-flag on $\pn{N}$ in good position
wrt to $\Sigma$. Then 

\n
($i$) We have the short exact sequence
in $\pe_Y^{\Sigma'_Y} \subseteq \pe_Y:$
$$
0 \lorw Q_{Y_{-n} -\emptyset } [-n] \lorw
\ldots \lorw  
Q_{Y_{-1}-Y_{-2} }[-1] \lorw Q_{Y-Y_{-1}} [0] \lorw Q \lorw 0;
$$
($i'$) If, in addition, $Y$ is affine, then
$ H^r(Y, Q_{Y_p -Y_{p-1}} ) =0, \;\; \forall r \neq p.
$

\n
($ii$) We have the short exact sequence
in $\pe_Y^{\Sigma'_Y} \subseteq \pe_Y:$
$$
0 \lorw 
Q 
\lorw  
{k_0}_* k_0^! Q  
\lorw 
 {k_{-1}}_* k_{-1}^! Q    [1]
 \lorw
 \ldots 
 \lorw 
 {k_{-n}}_* k_{-n}^! Q [n]
 \lorw 
 0;
$$
($ii'$) If, in addition, $Y$ is affine, then
$ H^r_c(Y, {k_{-p}}_* k_{-p}^! Q ) =0, \;\; \forall r \neq p.
$

\n
Let $({\frak F}, {\frak F}')$ be a pair of linear $n$-flags 
which are in good position with respect to $\Sigma$ and to each other.
Then

\n
($iii$) the  single complex $Gr^*_{\delta} Q [*]$  associated
with the double complex of perverse sheaves
$Gr^i_F Gr^j_G Q [i+j]$ is canonically isomorphic
to $Q$ in $D^b (\pe_Y).$ 

\n
($iii'$)  $H^{r} (Y,  ({k'_{-j}}_* {k'}_{-j}^! Q)_{Y_{i}-Y_{i-1} } )
= H^{r}_c (Y,  ({k'_{-j}}_* {k'}_{-j}^! Q)_{Y_{i}-Y_{i-1} } ) =0$
for every $r \neq i+j$.

\end{lm}
{\em Proof.}
Note that, in  Lemma \ref{brs}, we have $Y=Y_0$,
$j_!j^* Q = {k_0}_! k_0^* Q = Q_{Y_0-Y_{-1}}$, and
$j_* j^! Q = {k_0}_* k_0^! Q.$
More generally, we have $R\Gamma_{Z_{-j}-Z_{-j-1}} =
{k'_{-j}}_* {k'}_{-j}^!$ and $(-)_{Y_i - Y_{i-1}}=
{k_{i}}_! {k}_{i}^*.$

\n
Statement ($i$)  follows by a simple iteration of 
Lemma \ref{brs},
where one uses at each step
the fact that $\frak F$ is in good position wrt the initial
$\Sigma.$  In this step, the relative position
of the linear sections and the strata at infinity is unimportant.

\n
Statement ($i'$) follows from an iterated use
 of  Theorem \ref{swl} and Remark \ref{mm0}.
 Here it is important that the linear sections
 meet the strata at infinity transversally.

\n
Statements ($ii$) and ($ii'$) are proved in a similar way.

\n
The double complex is obtained as follows:
first resolve $Q$ as in ($ii$), then resolve
each resulting entry as in ($i$).
We thus have quasi isomorphism in $C^b (\pe_Y)$:
$Q \to Gr^*_G Q [*] \leftarrow Gr^{\bullet}_{\delta} Q
[\bullet]$ and ($iii$) follows.

\n
Finally, ($iii'$) now follows from (\ref{dova}).
\blacksquare

\begin{rmk}
\label{exfar}
{\rm
The formation of the left, right and bi-sided  resolutions
of $Q \in \pe_Y^{\Sigma_Y}$ in Lemma \ref{cellu} are  exact functor
with values in $C^b(\pe_Y).$
}\end{rmk}

\begin{ass}
\label{basassf}
{\rm  ({\bf Choice of the   pair of linear flags $\frak F, {\frak F}'$})
We fix a pair of  linear $n$-flags  $\frak F, {\frak F}'$ on $\pn{N}$ 
 in good position wrt $\Sigma$ and to each other.
A general pair in $F(N,n) \times F(N,n)$ will do.}
\end{ass}
\begin{rmk}\label{tuch}
{\rm Since 
$K \in D^{\Sigma_Y}_Y,$
the perverse sheaves
$\pc{b}{K} \in D_Y^{\Sigma_Y}$ and, with our choice
of the pair  $(\frak F, {\frak F}')$,
the conclusions of Lemma \ref{cellu} hold
for all the $\pc{b}{K}$.
}
\end{rmk}

\begin{lm}
\label{stk}
If $Y$ is affine, then 
$$
H^{r} ( Gr^a_F Gr^b_P  \, R\Gamma (Y,K)) =0, \qquad
\forall r \neq a-b,
$$
$$
H^{r} ( Gr^a_G Gr^b_P  \, R\Gamma_c (Y,K)) =0, \qquad
\forall r \neq a-b.
$$
If $Y$ is quasi projective, then
$$
H^r (Y, Gr^a_{\delta} Gr^b_P  R\Gamma (Y,K)) =
H^r (Y, Gr^a_{\delta} Gr^b_P  R\Gamma_c (Y,K)) =0
\qquad \forall r \neq a-b.
$$
\end{lm}
{\em Proof.}
We prove the first assertion, the second and third
 are proved in a similar way.
By (\ref{btfcm}),
 the group in question is
$$
H^{r- (a-b)} (Y, \pc{-b}{K}_{Y_a -Y_{a-1}}[a])
$$
and the required vanishing follows from 
Assumption \ref{basassf}, Remark \ref{tuch},   and Lemma \ref{cellu}.($i'$).
\blacksquare

\subsection{Proofs of  Theorems  \ref{tma}, \ref{tmac},  \ref{tmapl}
and \ref{cqpp}}
$\,$

\n
{\bf Proof of Theorems \ref{tma} and \ref{tmac}.}

\n 
By  the first two assertions of Lemma \ref{stk}, we can apply Proposition
\ref{bfpr} to $(R\Gamma (Y, K), P,F)$ and
to $(R\Gamma_c (Y, K), P,G)$   and conclude.
\blacksquare

\bigskip
\n
{\bf Proof of Theorem \ref{tmapl}.}

\n
We prove the version for $H(X,C)$. The case of $H_c(X,C)$
is proved in a similar way.

\n
Given the  fixed embedding $Y \subseteq \pn{N}$,
pick a
a stratification $\Sigma$ of $\pn{N}$ adapted to
$f_*C$ and to the embedding. Choose a linear $n$-flag
$\frak F$ on $\pn{N}$ in good position wrt $\Sigma$,
e.g. general.
Let $Y_*$ be the corresponding $n$-flag on $Y$ and 
set $X_*:= f^{-1}Y_*.$ Denote by $\widetilde{i},\widetilde{j},
\widetilde{k}$ 
the associated embeddings as in (\ref{ijk}).

\n
By Theorem \ref{tma},  
the perverse spectral sequence for $H(Y, f_*C)$, i.e. the perverse
Leray spectral sequence for $H(X,C)$, is the shifted $Y_*$ flag
spectral sequence for $H(Y,f_*C)$  ($F$-version).

\n
Our goal is to identify the
$Y_*$ flag spectral sequence for $H(Y,f_*C)$ with the
$X_*$ spectral sequence for $H(X,C)$.
It is sufficient to show that
\begin{equation}
\label{357}
(R\Gamma (X, C), F_{X_*} ) = (R\Gamma (Y, f_* C), F_{Y_*});
\end{equation}
in fact,  the shifted versions would  also coincide
and we would be done.
 In general,  the two filtered complexes  
 for $Y_*$ and $X_*$ do not coincide,
due to the failure of the base change theorem. In the present case,
transversality prevents this from happening.

\n
We assume that $C$ is injective. The filtered complex
$(C, F)$ is of $c$-soft type.   
On   varieties
$c$-soft and soft are equivalent notions and  
soft sheaves are $f_!$ and $f_*$-injective.

\n
We have the filtered complex
 $(R \Gamma (X, C), F_{X_*})$, i.e. the result of applying
 $\Gamma (X,-)$ to the $C$-analogue of (\ref{ttff}).

\n
Transversality ensures that we have the first equality
in (\ref{bctcms}):
$ {\widetilde{j_p}}_! \widetilde{j}_p^* f_* C =f_* {j_p}_! j_p^* C.$
This implies that, by applying $f_*$ to the $C$-analogue of (\ref{ttff}),
we obtain the $f_*C$ analogue of (\ref{ttff}) on $Y$ wrt to $Y_*$, 
i.e. (\ref{357}) holds and we are done.\blacksquare

\bigskip
\n
{\bf Proof of Theorem \ref{cqpp}.}

\n
In view of  the third assertion of Lemma \ref{stk},
the proof is analogous to the proofs given above. 
\blacksquare

\subsection{Resolutions in $D^b(\pe_Y)$}
\label{resindbp}
In an earlier version of this paper, we  worked
in the derived category of perverse sheaves $D^b(\pe_Y)$
 which, in the case of field coefficients, is equivalent to $D_Y$ 
 (\ci{be}).
  We are very thankful to one of the anonymous referees
for suggesting the  considerably more elementary 
approach contained in this paper which takes place in $D(Ab)$.  
On the other hand, the approach in $D^b(\pe_Y)$ 
explains the relation $``P=Dec(F)"$ at the level of complexes
of (perverse) sheaves, i.e. before taking cohomology.
 We outline this approach in the
case of the $F$-construction on  $Y$ affine.  
We omit writing down the similar
details in the case of the $G$-construction in the affine 
case and in the case  of the $\delta (F,G)$-construction
in the quasi projective case. 

The approach is based on Beilinson's Equivalence
Theorem \ci{be}.

In what follows,  $Y$ is affine, we work with field
coefficients, e.g. $\rat$,  $D^b(\pe_Y)$ is endowed with
the standard $t$-structure, $D_Y$ with the perverse $t$-structure.
An equivalence of $t$-categories is a
functor between triangulated categories with $t$-structures which is
additive, commutes with translations, preserves distinguished
triangles, is $t$-exact (i.e. it preserves the hearts) and is an equivalence.

\begin{tm}
\label{beitm} {\rm (\ci{be})}
There is an equivalence of $t$-categories, called
the realization functor
$$
r_Y:  D^b(\pe_Y) \stackrel{\simeq}\lorw D_Y.
$$
\end{tm}

An outcome of this result is that it implies that, up to replacing
$K \in D_Y$ with a complex naturally isomorphic to it,
 there is a filtration $B$ on $K$ such that
  $Gr_B^b K[b] \in \pe_Y$. Recalling
   $\S$\ref{canlift}, this means that
 $(K,B)$  is in
 the heart   $D_Y F_{\beta}$ of the canonical lift
  to $D_YF$
  of the perverse $t$-structure on $D_Y$.
 This circumstance, coupled with the construction
(\ref{cfses}), allows to describe an inverse $s_Y$ to $r_Y$,
i.e. to assign to $K$ a  complex of perverse sheaves
$$
s_Y(K) = s_Y (K,B)= Gr_B^* K [*] = : {\cal K}^*  \in D^b(\pe_Y).
$$

Fix a stratification $\frak S$
  of $Y$ such that all the finitely many non-zero $Gr^b_BK$
  are $\frak S$-constructible. 
  Note that if $K \in D^{\Sigma'}_Y$, then it is possible
  that 
  $Gr^b_BK \notin D^{\Sigma'}$ and one may need to refine.
   Choose an embedding $Y \subseteq \pn{N}$, 
   a stratification $\Sigma $
  on $\pn{N}$ adapted (cf.
  Definition \ref{adpted}) to  $\frak S$ and to  the embedding,
  and a linear $n$-flag $\frak F$ on $\pn{N}$ in good position
  (cf. Definition \ref{basdef}; a general one will do)
  with respect to $\Sigma$. 
  
  Let  $\Delta = \Delta (F,P)$ be the diagonal filtration.
  By transversality, we have that $(K,\Delta) \in D_Y F_{\beta}.$
  We obtain the double complex
  ${\cal K}^{*,*}:= Gr^*_F Gr^*_P K$, 
   with associated single complex $s({\cal K}^{*,*})^*
   =s_Y( K, \Delta)$ that maps
    quasi isomorphically onto ${\cal K}^*.$
   We also have $H^r (Y, {\cal K}^{*,*}) =0$ for every
   $r \neq 0$ so that  we have  obtained 
   a resolution with    $H(Y,-)$-acyclic   entries.

 The single complex $s({\cal K}^{*,*})$ admits the
 b\^{e}te filtration by rows $B_{row}$, where $B_{row}^q
 s ( {\cal K}^{*,*} )$ is the single complex associated with the 
 double complex ${\cal K}^{*, *\leq q}$.
 i.e. the result of replacing with zeroes the entries
 strictly above the $q$-th row.  
 
 There is another filtration, $Std_{col},$ where 
 $Std_{col}^p  \,  s ( {\cal K}^{*,*} )$ is the single complex associated
 with the double complex obtained by keeping the columns $p' <-p,$
 replacing the columns $p'>-p$ with zeroes,
 and replacing the entries ${\cal K}^{p,q}$ in  $(-p)$-th column 
 by $\ke \, \{ {\cal K}^{p,q} \to {\cal K}^{p+1, q} \}.$ 
 
By Remark \ref{exfar},
the exactness properties of the construction of the resolution
 of Lemma \ref{cellu}  ensure that
 the natural map $( s ( {\cal K}^{*,*} ), Std_{col}) \to ({\cal K}^*, Std)$
 is a filtered quasi isomorphism.
 
It is via  this construction
that the relation $Dec(F) = P$ becomes transparent: 
it holds in  $D^b(\pe_Y) \simeq D_Y$ and it descends to
$D(Ab)$: 

 \n
 1) it is elementary to verify that $Dec(B_{row}) = Std_{col}$
 (cf. \ci{beabs}, Remark 3.11.1);

 \n
2) 
the filtered complex of perverse sheaves
  $ (s ( {\cal K}^{*,*} ), B_{row}) $ corresponds to
  $(K,F)$ under the  equivalence $r_Y$;
  
  \n
 3) 
 the $t$-exactness of  $r_Y$ ensures that
  the filtered complex  $ (s ( {\cal K}^{*,*} ), Std_{col}) \simeq
 ({\cal K}, Std) $ correspond to $(K,P)$;
  
  \n
  4) by the exactness of the construction,
 the complex $s( H^0(Y, {\cal K}^{*,*}))$ inherits
 the relation
  $Dec(B_{row}) = Std_{col}$.

  \n
  5) 
  the bifiltered complex 
  $(s( H^0(Y, {\cal K}^{*,*})),  Std_{col}, B_{row})$
  realizes $(R\Gamma (Y, K), P, F)$
  and 4) and  
   the perverse spectral sequence is identified
  with the shifted flag spectral sequence.

 \begin{rmk}
 \label{lderrder}
 {\rm
 On the affine $Y$,  the functor $H^0_{\pe_Y}: \pe_Y \to Ab$,
 $Q \mapsto H^0(Y,Q)$ is  right-exact.
 By \ci{be}, $\S$3, this right-exact functor admits
 a left-derived functor  $LH^0_{\pe_Y}: D^-(\pe_Y) \to
 D^-(Ab)$. The complex in 
 Step 4)   realizes $LH^0_{\pe_Y} ({\cal K})$. 
 }
 \end{rmk}

\section{Applications}
\label{appil}
The following results are due to M. Saito \ci{msaito}
who used his own mixed Hodge modules.
We offer a proof based on the methods of this paper. 
\begin{tm}
\label{ap1}
Let $Y$ be quasi projective.
The perverse spectral sequences
$$
E_1^{p,q} = H^{2p+q}(Y, \pc{-p}{\zed_Y}) \Longrightarrow
H^*(Y, \zed),
$$
$$
E_1^{p,q} = H^{2p+q}_c(Y, \pc{-p}{\zed_Y}) \Longrightarrow
H^*_c(Y, \zed)
$$
are spectral sequences in the category of mixed Hodge structures.
\end{tm} 
{\em Proof.}
We prove the first statement when $Y$ is affine.
 The other variants are proved in  similar
ways. 
By Theorem \ref{tma}, there is a $n$-flag $Y_*$ on $Y$
such that the perverse spectral sequence for $H(Y, \zed)$
is  the shifted   spectral sequences of the 
flag spectral sequence
$$
E_1^{p,q} = H^{p+q}(Y_p, Y_{p-1},\zed) \Longrightarrow
H^*(Y, \zed),
$$
which is in the category of mixed Hodge structures.
\blacksquare

\begin{tm}
\label{ap12}
$f:X \to Y$ be a map of varieties with $Y$ quasi projective.
The perverse Leray spectral sequences
$$
E_1^{p,q} = H^{2p+q}(Y, \pc{-p}{f_*\zed_Y}) \Longrightarrow
H^*(X, \zed),
$$
$$
E_1^{p,q} = H^{2p+q}_c(Y, \pc{-p}{f_!\zed_Y}) \Longrightarrow
H^*_c(X, \zed)
$$
are spectral sequences in the category of mixed Hodge structures.
\end{tm} 
{\em Proof.}
We prove the case of cohomology over an affine 
base $Y$ and leave the rest to the reader.
By Theorem \ref{tmapl}, 
the perverse Leray spectral sequence for $H(X,\zed)$
is the  shifted $X_*$ flag spectral sequence wrt a suitable 
$n$-flag $X_*$ on $X$.
This latter is in the category of mixed Hodge structures.
\blacksquare

\begin{rmk}
\label{marenostrum}{\rm 
({\bf Mixed Hodge structures and the decomposition theorem})
In the paper \ci{decmightam}, we endow  the 
cohomology of the direct summands
appearing in the decomposition theorem 
for the proper push forward of the intersection cohomology complex
of a proper variety, with natural pure polarized Hodge structures.
These structures 
arise as subquotients of the pure Hodge structure
of the cohomology of a resolution 
of the singularities of the domain of the map.
In particular, this endows the intersection cohomology groups
of proper varieties with pure polarized Hodge structures.
In the paper \ci{decmigseattle}, we prove that for projective morphisms
of projective varieties,
one can realize the direct sum splitting mentioned above
 in the category of pure Hodge structures.
The methods of this paper allow to  endow 
the intersection cohomology groups
$I\!H(Y,\zed)$ and $I\!H_c(Y,\zed)$ of a quasi projective variety
with a mixed Hodge structure and to extend all the results
of \ci{decmightam} to the case of quasi projective varieties.
Furthermore, we  compare the resulting mixed Hodge structures
with the ones arising from M. Saito's work and we show that they coincide.
Details will appear  in \ci{pfII}.
}\end{rmk}

Authors' addresses:

\smallskip
\n
Mark Andrea A. de Cataldo,
Department of Mathematics,
 Stony Brook University,
Stony Brook,  NY 11794, USA. \quad 
e-mail: {\em mde@math.sunysb.edu}

\smallskip
\n
Luca Migliorini,
Dipartimento di Matematica, Universit\`a di Bologna,
Piazza di Porta S. Donato 5,
40126 Bologna,  ITALY. \quad
e-mail: {\em migliori@dm.unibo.it}

\end{document}